\documentclass[journal]{IEEEtran}
\usepackage{amsmath,amsfonts,amssymb}
\usepackage{cite}
\usepackage{color}
\usepackage{graphicx}
\usepackage{array}

\begin{document}

\title{Control Strategy Design for Power Quality Management in 
Active Distribution Networks}

\author{%
Surya~Chandan~Dhulipala,~\IEEEmembership{Student Member,~IEEE,}
Arturo~Bretas,~\IEEEmembership{Senior Member,~IEEE,}
and Enrique~Baeyens,~\IEEEmembership{Member,~IEEE,}
\thanks{This work was supported in part by the University of Florida
and NSF EECS-1129061}
\thanks{S.C. Dhulipala, A. Bretas, S. Meyn and P.P. Khargonekhar
are with Electrical and Computer Engineering Department, 
University of Florida, FL, USA.}%
\thanks{E. Baeyens is with
Instituto de las Tecnolog\'{\i}as Avanzadas de la Producci\'on,
Universidad de Valladolid, Valladolid, Spain}%
\thanks{P. Barooah is with Mechanical and Aerospace Engineering
Department, University of Florida, FL, USA}%
}

\maketitle

\begin{abstract}
Uncertainty and variability of renewable energy sources
present an imperative technical challenge for Electrical Distribution
Utilities. Power Quality (PQ) indices represent quality of energy
delivered and reliability. In this work, a control strategy design for
real-time PQ management in active distribution systems is presented. The
present work addresses management of voltage
fluctuations induced by the variability in renewable generation. A ``zero
energy reserve'' approach to tackle fluctuations of renewable energy
generators is developed. The power consumption of flexible loads is
modulated to reduce the technical losses and peak load of the feeder in
this reactive power (VAR) control strategy. A Volt/VAR control strategy
formulation, using the capability of smart inverters to provide dynamic
reactive power is presented.  Unbalanced power flow, different load
profiles and flexible loads as ``virtual energy storages'' are used to
improve voltage profile and reduce technical losses while maintaining
system reliability. IEEE 13 bus distribution system is used for control
strategy design validation. Comparative results indicate reduction in the
system technical losses and the stress on automatic voltage
regulators (AVR). The ease of design of control strategy indicate
potential real-life application.
\end{abstract}

\begin{IEEEkeywords}
Distributed power generation, 
Energy storage,
Markov processes, 
Reactive power control.
\end{IEEEkeywords}

\section{Introduction}
\label{sec:1}

Voltage regulation (VR) in distribution systems is traditionally
performed at distribution network substation level. Capacitor banks
adapt to the reactive power requirement of the system on a slower
timescale (typically every hour). Due to uncertainty of 
available power from renewable energy sources and the upward trend of
solar plant penetration into the distribution grid, there is a need for
VAR control on a faster time scale \cite{moreno2007}. 
Smart inverters, interfacing the
solar plants at the point of interconnection are excellent controllable
resources, which can follow a remote command to change their power factor.
Three-phase inverters are ultimate reactive power generators, with
response times of $1$~ms to full output \cite{schauder2014}. 
Flexible loads in this document are defined as loads that can follow 
a regulation signal, providing ancillary services to the network. 
In \cite{lin2015}, the authors show that HVAC systems in commercial 
buildings can provide demand side system frequency regulation.

Comprehensive studies have been made on Volt/VAR control for
distribution systems. Most works consider the problems of optimal
sizing, placement and switching schedules of capacitor banks and on-load
tap changers (OLTC) \cite{baran1989}, \cite{liu2009} and \cite{liu2002}. 
Inverter control on a fast timescale has been studied in some recent works 
\cite{farivar2011}, \cite{yeh2012} and \cite{farivar2012}.
In \cite{farivar2011},
authors present a Volt/VAR control problem that minimizes both technical
losses and power consumption. The core idea of this strategy was
minimization of the weighted sum of voltages along the feeder 
using a balanced power flow model. In distribution
systems, the unbalanced power flow is a model closer to reality 
than the DistFlow equations introduced in \cite{baran1989} and
used in \cite{farivar2011}. 
The potential benefits of flexible loads and real time
power quality management were also not addressed in \cite{farivar2011}.

In this work, a centralized control strategy for real-time PQ
management on distribution networks is presented. The present control 
strategy optimizes the technical losses, deviation from nominal voltage 
and power consumption of flexible load in a cost-effective approach
considering the inherent distribution network characteristics.
The inherent unbalanced operation, different
load types, Markov chain modelling for solar generation and load
forecasting, modulation of power consumption of flexible loads and the
issues of following the fluctuations of renewable energy generation are
considered in this work. The traditional VAR control using capacitor
banks on a slow timescale has been studied in previous works. In this work,
the control of smart inverters on a short timescale for real-time PQ
management is investigated. 

The remainder of the paper is organized as follows. In Section~\ref{sec:2}, 
we present the problem formulation. In Section~\ref{sec:3}, the fundamental
aspects of Markov chain modelling of solar energy and load forecasting
are addressed. In Section~\ref{sec:4}, a zero energy reserve approach is
proposed and the control of an Energy Storage System (ESS) to follow the
output fluctuations of renewable energy generators is discussed. In
Section~{\ref{sec:5}, the proposed power quality control strategy design is
reported. Comparative simulation results, discussion and analysis are
presented in Section~\ref{sec:6}. Finally, some conclusions are given in
Section~\ref{sec:7}.

\section{Problem Formulation}
\label{sec:2}

\subsection{Control Problem Statement}

Volt/VAR regulation is a control problem on two timescales. Slow
timescale control is done hourly (every $t$) and inverter fast timescale control
is done every minute (every $t'$). Current SCADA limits the timescale of control 
strategy to a minute. In this work, the inverter power factor which can be
changed dynamically to match the uncertainty of renewable energy
generators will be used for fast timescale Volt/VAR control. The
configuration of capacitor banks is changed in the slow timescale
control. This will reduce the number of tap changes of OLTC, thereby
increasing their lifecycle.

\subsection{Distribution Grid Model}

Consider a distribution network with $n$ buses, where $\mathcal B$ is a 
set of all the buses on the system $\mathcal B = (1,2,\ldots,n_B)$. Let
$\Phi=\{a,b,c\}$ the set of phases.
Consider the distribution network that utilizes capacitor
banks and tap changing transformers for voltage regulation. The state of
the capacitor banks and tap changing transformers at time $t$ is modeled
by an integer vector $s(t)\in \mathcal A$, where $\mathcal A$ is the set of feasible
values of vector $s(t)$. The admittance matrix of the grid at time $t$ is
given by $Y(s(t))$. The elements of the admittance matrix
$Y_{ij-km} = G_{ij-km} + j B_{ij-km}$ are
\begin{align}
\label{eq:1}
Y_{ij-kk} &=
bsh_{i-k} + \sum_{m\in\Omega_k} 
\left( y_{ij-km} + bsh_{ij-km} \right) \\
\label{eq:2}
Y_{ij-km} &= 
- y_{ij-km}
\end{align}
where $\Omega_k$ is the set of buses adjacent to bus $k$,
$y_{ij-km}$, $bsh_{ij-km}$ are series admittance and susceptance,
respectively, of lines connected between buses $(i,j) \in \mathcal B$ 
and phases $(k,m) \in \Phi$,
and $bsh_{i-k}$ is the shunt capacitance at bus $i$, phase $k$ which
when there is presence of capacitor bank, is determined by state
$s(t)$ of the capacitor banks.

\subsection{Generator Model}

The complex power injected by the aggregated generators at node $k$ and
time $t$ is denoted as $S_k^G(t) = P_k^G(t)+jQ_k^G(t)$. 
The active and reactive power that each generator can provide is bounded by 
\begin{align}
&P_k^{G,\min}(t) \leq P_k^G(t) \leq P_k^{G,\max}(t) \\
&Q_k^{G,\min}(t) \leq Q_k^G(t) \leq Q_k^{G,\max}(t)
\end{align}
If there is no generator at bus $k$, then $P_k^G(t) = Q_k^G(t)  = 0$. 
Time dependence of the active power takes into account generation from
renewable sources that is time-varying. If the generator at bus $k$ uses
only renewable energy sources, it is assumed in this work that the
generation profile $P_k^G(t)$ is known/estimated in advance for the interval
of time under study. Also, since one cannot schedule renewable
resources, $P_k^{G,\min}(t) = P_k^{G,\max}(t)$ is assumed.

\subsection{Unbalanced Power Flow Formulation}

Let $P^{\exp}$ and $Q^{\exp}$ denote the expected/estimated net injected 
real and reactive power while $P$ and $Q$
denote the realized injected real and reactive power. The difference between
the estimated and the realized real and reactive powers is given by the
following equation:
\begin{align}
\begin{bmatrix}
\Delta P \\ \Delta Q
\end{bmatrix} &=
\begin{bmatrix}
P^{\mathrm{exp}} \\ Q^{\mathrm{exp}}
\end{bmatrix} -
\begin{bmatrix}
P \\ Q
\end{bmatrix}
\label{eq:5}
\end{align}

The goal is to estimate system states (bus voltage magnitudes and
angles) which minimizes functions $\Delta P$ and $\Delta Q$ to a 
pre-established convergence tolerance. 
Realized power in (\ref{eq:5}) can be given by the following 
expressions:
\begin{align}
P_{i-k} &= 
V_{i-k} \sum_{j\in\mathcal B} \sum_{m\in \Phi}
V_{j-m} (
G_{ij-km} \cos \theta_{ij-km} + \nonumber \\
& \quad B_{ij-km} \sin \theta_{ij-km}
) \label{eq:6} \\
Q_{i-k} &= 
V_{i-k} \sum_{j\in\mathcal B} \sum_{m\in \Phi}
V_{j-m} (
G_{ij-km} \sin \theta_{ij-km} - \nonumber \\
& \quad B_{ij-km} \cos \theta_{ij-km}
) \label{eq:7}
\end{align}
where $P_{i-k}$ and $Q_{i-k}$ represent the active and reactive power
injections at bus $i$, phase $k$.
$V_{i-k}$ and $V_{j-m}$ represent the magnitude of voltage phasor at bus $i$, 
phase $k$ and bus $j$, phase $m$, respectively.
$G_{ij-km}$ is the conductance of the lines connected between bus $i-j$ 
and phase $k-m$, respectively.
$B_{ij-km}$ is given as in (\ref{eq:1}) and (\ref{eq:2}).
$\theta_{ij-km}$ is the angular difference between 
bus $i$, phase $k$ and bus $j$, phase $m$.

The system of non-linear algebraic equations modelling the power flow
in the distribution grid  is solved using Newton-Raphson (NR) approach. 
The initial system state is normally considered as a flat voltage profile. NR
approach solves (\ref{eq:5}) iteratively, by linearizing the system 
at each time step. The linearized system of equations is given by

\begin{align}
\label{eq:8}
J
\begin{bmatrix}
\Delta \theta \\ \Delta V 
\end{bmatrix}
&=
\begin{bmatrix}
\Delta \theta \\ \Delta V 
\end{bmatrix}
\end{align}
with the states being updated according to the following equation,

\begin{align}
\label{eq:9}
\begin{bmatrix}
\theta \\ V 
\end{bmatrix}^{\nu+1} &=
\begin{bmatrix}
\theta \\ V 
\end{bmatrix}^{\nu} +
\begin{bmatrix}
\Delta \theta \\ \Delta V 
\end{bmatrix}^{\nu}
\end{align}

The process is repeated until functions 
$\Delta P$ and $\Delta Q$ are smaller than a pre-established convergence 
tolerance. 

The Jacobian matrix can be derived by considering different types
of loads. Consider the load model given by the following equations
\begin{align}
\label{eq:10}
P_{k-i}^{\mathrm{exp}} &=
P_{k-i}^{\mathrm{nom}}
\sum_{l=1}^n a_l 
\left( 
\frac{V_{k-i}}{V_{k-i}^{\mathrm{nom}}}
\right)^{{n_p}_l} \\
\label{eq:11}
Q_{k-i}^{\mathrm{exp}} &=
Q_{k-i}^{\mathrm{nom}}
\sum_{l=1}^n b_l 
\left( 
\frac{V_{k-i}}{V_{k-i}^{\mathrm{nom}}}
\right)^{{n_q}_l}
\end{align}
where $P_{k-i}^{\mathrm{exp}}$ and $Q_{k-i}^{\mathrm{exp}}$ are
expected real and reactive power of bus $k$ for a given
voltage $V_k$, $V_k^{\mathrm{nom}}$ is the bus $k$ nominal voltage, 
$a_l$ and $b_l$ are the bus $k$ expected percentages of  
$\left(P_{k-i}^{\mathrm{exp}},Q_{k-i}^{\mathrm{exp}}\right)$
with respect to coeficients $n_{p_l}$ and $n_{q_l}$, with
$\sum_{l=1}^{\ell} a_l =1$, $\sum_{l=1}^{\ell} b_l = 1$, 
where $\ell$ denotes the different load types present of the bus.

The Jacobian matrix is:
\begin{align}
\label{eq:12}
J &=
\left[
\begin{array}{cc}
\displaystyle
\frac{\partial P}{\partial \theta} &
\displaystyle
\frac{\partial P}{\partial V} \\[2ex]
\displaystyle
\frac{\partial Q}{\partial \theta} &
\displaystyle
\frac{\partial Q}{\partial V}
\end{array}
\right]
\end{align}
and its elements are:

\begin{align}
\frac{\partial P_{i-k}}{\partial \theta_{i-k}} &= 
-V_{i-k}^2 B_{ii-kk}-Q_{i-k}
\label{eq:13}
\end{align}

\begin{align}
\frac{\partial P_{i-k}}{\partial \theta_{j-m}} &= 
V_{i-k}V_{j-m} ( G_{ij-km}\sin \theta_{ij-km} - \nonumber \\
& \quad
B_{ij-km} \cos \theta_{ij-km} )
\label{eq:14}
\end{align}

\begin{align}
\frac{\partial P_{i-k}}{\partial V_{j-m}} &= 
V_{i-k} ( G_{ij-km}\cos \theta_{ij-km} + \nonumber \\
& \quad
B_{ij-km} \sin \theta_{ij-km} )
\label{eq:15}
\end{align}

\begin{align}
\frac{\partial Q_{i-k}}{\partial \theta_{i-k}} &= 
-V_{i-k}^2 G_{ii-kk} + P_{i-k}
\label{eq:16}
\end{align}

\begin{align}
\frac{\partial Q_{i-k}}{\partial \theta_{j-m}} &= 
-V_{i-k}V_{j-m} ( G_{ij-km}\cos \theta_{ij-km} + \nonumber \\
& \quad
B_{ij-km} \sin \theta_{ij-km} )
\label{eq:17}
\end{align}

\begin{align}
\frac{\partial Q_{i-k}}{\partial V_{j-m}} &= 
V_{i-k} ( G_{ij-km}\sin \theta_{ij-km} - \nonumber \\
& \quad
B_{ij-km} \cos \theta_{ij-km} )
\label{eq:18}
\end{align}

\begin{align}
\frac{\partial P_{i-k}}{\partial V_{i-k}} &= 
V_{i-k} G_{ii-kk} + \frac{P_{i-k}}{V_{i-k}} - \nonumber \\
& \quad
P_{k-i}^{\mathrm{nom}} \sum_{l=1}^n {n_p}_l a_l
\left(
\frac{V_{k-i}}{V_{k-i}^{\mathrm{nom}}}
\right)^{{n_p}_l-1}
\label{eq:19}
\end{align}

\begin{align}
\frac{\partial Q_{i-k}}{\partial V_{i-k}} &= 
-V_{i-k} B_{ii-kk} + \frac{Q_{i-k}}{V_{i-k}} - \nonumber \\
& \quad
Q_{k-i}^{\mathrm{nom}} \sum_{l=1}^n {n_q}_l b_l
\left(
\frac{V_{k-i}}{V_{k-i}^{\mathrm{nom}}}
\right)^{{n_q}_l-1}
\label{eq:20}
\end{align}

\section{Solar Forecasting and Load Modeling: A Markov Chain Approach}
\label{sec:3}

\subsection{Introduction}

A Markov chain is a probabilistic model that represents system dynamics
with a finite number of states \cite{Norris1998}, \cite{Bremaud1999}. 
A Markov chain is characterized by a
pair of elements $(\mathcal N, \Pi)$  where  
$\mathcal N = \{1,2,\ldots,N\}$ is the set of finite states and 
$\Pi$ is the matrix of transition probabilities. 
Let $X_t$ be the state of the Markov
chain at time $t$. The state can change according to the corresponding
transition probability $\pi_{ij}= \mathrm{Pr} \{X_{t+1}=j|X_{t}=i\}$. 
Thus, the matrix $\Pi$ has nonnegative entries and
the sum of each row is equal to one. 

\subsection{Solar Insolation}

\subsubsection{Markov chain modelling}

Markov chains have been previously used to model and simulate 
solar insolation \cite{ehnberg2005}, \cite{bright2015}. 
In our Markov chain model, 
using recorded measurements of solar insolation from NREL National
Radiation Data Base \cite{wilcox2008}, the maximum 
and minimum solar insolation profile for each
hour of each month of the year can be estimated. The variation range is
split in $N$ levels that are labeled by the elements of the state set
$\mathcal N$.
The upper level of radiation is labeled by $1$ and the lower level by $N$
and the transitions probabilities can be estimated from the data.
Therefore, the transition probabilities encode the system dynamic
characteristics. In this work, a total of twelve Markov chain models,
one for each month of the year, are used.

\subsubsection{Solar insolation simulation}

A simulation of the solar insolation for a given day is realized as
follows. Select the Markov chain model for the corresponding month of
the year and the variation range for each hour. Select the initial state
and run a simulation for each hour of the day. At each time instant $t$ a
random uniform noise is added such that the solar insolation is random
but remains in the corresponding level given by the Markov chain state.

\subsubsection{Solar insolation forecasting}

The Markov chain model can be applied to forecast the solar insolation.
A Markov chain corresponding to the month of the day that one wants to
forecast and the variation range for the corresponding hour is
considered. The forecasted solar insolation
$\hat P_{t+k}^S$ at time $t+k$ assuming that
the state at time $t$ is $X^S_t=i$ can be obtained using the conditioned
expected value of the future state $\mathbb E[ X^S_{t+k} | X^S_{t}=i]$
and the range of variation $R^S(t+k)$ for the hour of the day corresponding
to time instant $t+k$, \emph{i.e.}
\begin{align*}
\hat P^S_{t+k} &=
\frac{1}{N} R^S(t+k) \left( \mathbb E[ X^S_{t+k} | X^S_{t}=i] -
\frac{1}{2} \right) \\
&= \frac{1}{N} R^S(t+k) \left( \sum_{j=1}^N b_j^{\top} \Pi_S^k b_i -
\frac{1}{2} \right)
\end{align*}
where $b_i$ is the $i$-th column vector of the identity matrix $I_{N_S}$
of dimension $N_S$, and $\Pi_S^k = \Pi_S^1\cdot \Pi_S^{k-1}$ and
$\Pi_S^0=I_{N_S}$.

\subsection{Load Forecasting}

\subsubsection{Load Modelling}

The modeling, simulation and forecast of load variation can be done
similarly using a discrete Markov chain. 
Using data from the past, one can obtain the maximum and minimum load
for a certain time interval (typically one hour) of a given week day,
for example a Tuesday of June. These maximum and minimum values define
the variation range of the load for this hour. This range is split in $N$
levels and a state of the Markov chain is assigned to each level. The
state corresponding to the upper subinterval of load variation is
labeled by $1$ and the lower one by $N$. The transition probabilities are
estimated from the data.

\subsubsection{Simulating the loads}

A simulation of the load variation for a given day is realized as
follows. Select the Markov chain model and the variation range for each
time interval of the corresponding week day and month of the year.
Select an initial state and run a simulation for each hour of that day.
At each time instant the load is given by the realization of a random
variable whose expectation is the mean value of the subinterval
corresponding to the Markov chain state and with range or variation
given by the range of the subinterval

\subsubsection{Flexible Load Forecasting}

The Markov chain model can also be used to forecast the 
flexible load variation. For instance, the Markov chain corresponding to
the week day and the month that  one wants to forecast and the variation
range for the corresponding time can be selected for the corresponding
time interval. The forecasted load $\hat P^L_{t+k}$ at time $t+k$
assuming that the state at time $t$ is $X_t^L=i$ is obtained using
the conditioned expected value of the future state
$\mathbb E[ X^L_{t+k} | X^L_{t}=i]$ and
the range of variation $R^L(t+k)$ for the hour of the day corresponding
to
time instant $t+k$, \emph{i.e.}
\begin{align*}
\hat P^L_{t+k} &=
\frac{1}{N} R(t+k) \left( \mathbb E[ X_{t+k}^L | X_{t}^L=i] -
\frac{1}{2} \right) \\
&= \frac{1}{N} R^L(t+k) \left( \sum_{j=1}^N b_j^{\top} \Pi_L^k b_i -
\frac{1}{2} \right)
\end{align*}
where $b_i$ is the $i$-th column vector of the identity matrix $I_{N_L}$
of dimension $N_L$, and $\Pi_L^k = \Pi_L^1\cdot \Pi_L^{k-1}$ and
$\Pi_L^0=I_{N_L}$.

The load profiles of a $13.8$~KV feeder in Florida, which were used in
the simulations are illustrated in Fig.~\ref{fig:1} and 
Fig.~\ref{fig:2}. Different
load profiles were used to simulate different seasons (Fig.~\ref{fig:3}). 
The profiles of different loads were derived from the feeder load profiles 
using their generic models (obtained by using averaged load profile for
different types of loads).

At any given state of the load, the amount of load which can be modulated 
is determined by the minimum and maximum envelopes of load consumption 
\cite{lin2015} as shown in the Fig.~\ref{fig:2},
where $P_{\min}^L=\min R^L(t)$, $P_{\max}^L=\max R^L(t)$, and $R^L(t)$
is the variation range of the load at time $t$.

\begin{figure}
\centering
\includegraphics[width=1\columnwidth]{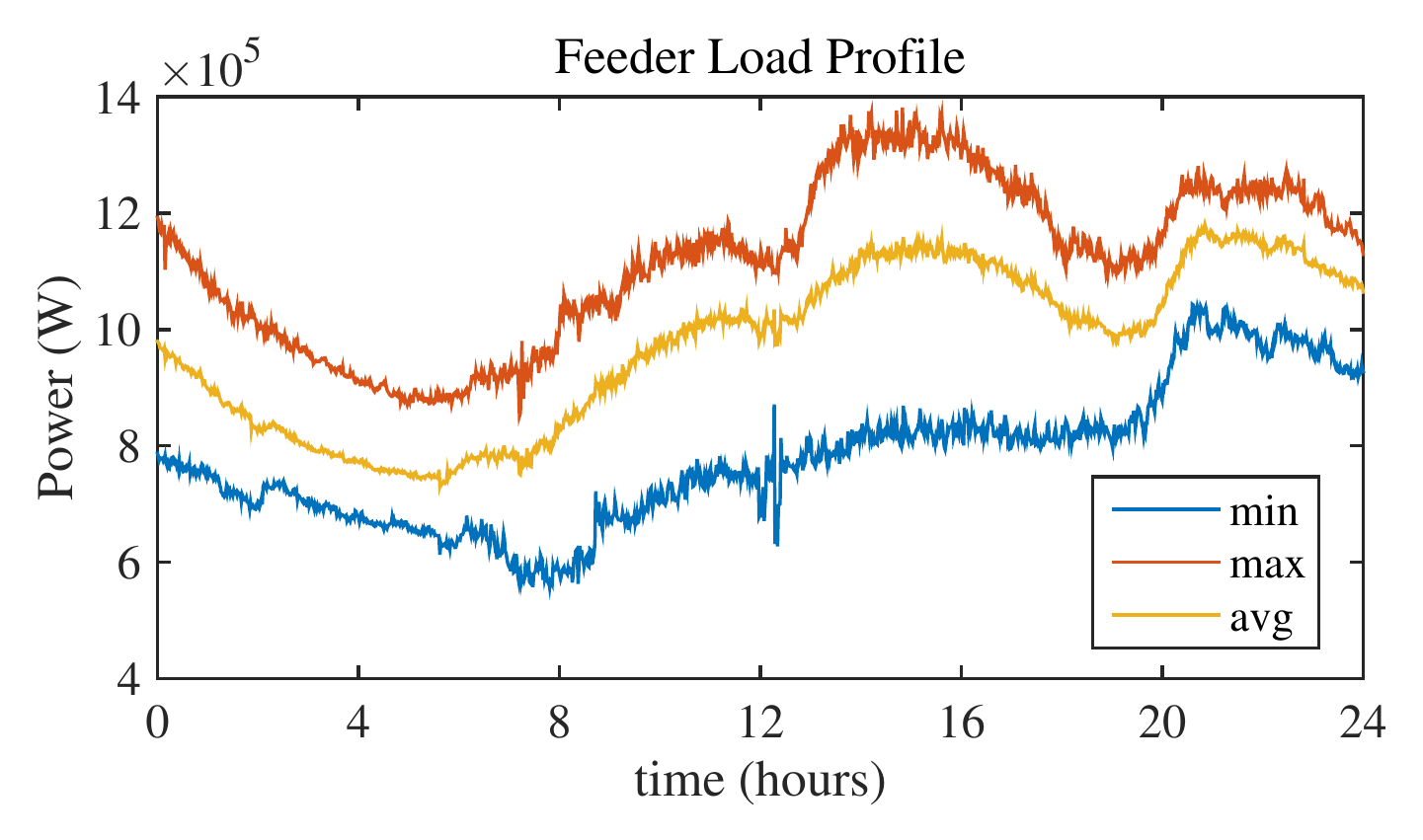}
\caption{Feeder Load Profile (Summer)}
\label{fig:1}
\end{figure}

\begin{figure}
\centering
\includegraphics[width=1\columnwidth]{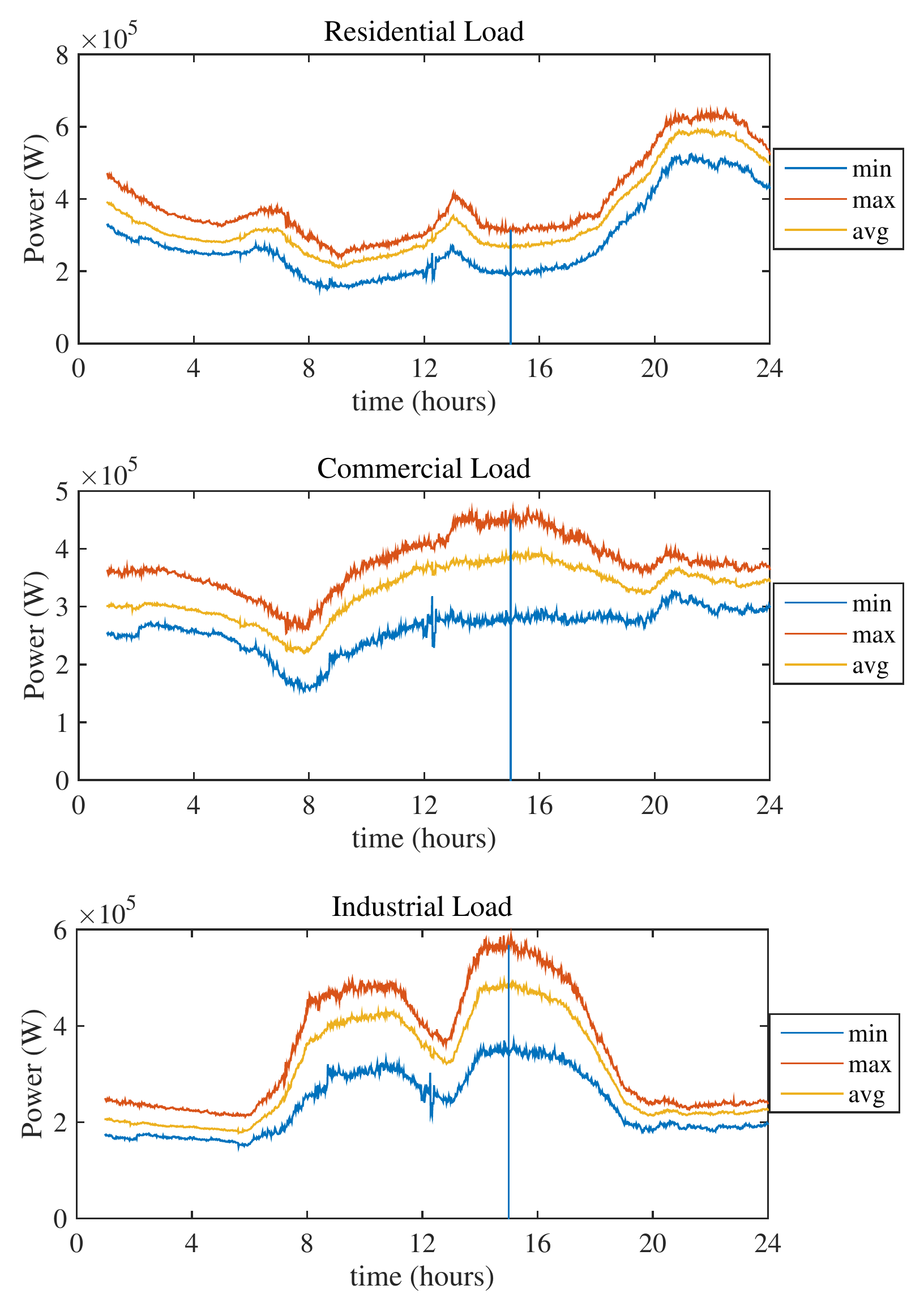}
\caption{Different Load Profiles (Summer)}
\label{fig:2}
\end{figure}

\begin{figure}
\centering
\includegraphics[width=1\columnwidth]{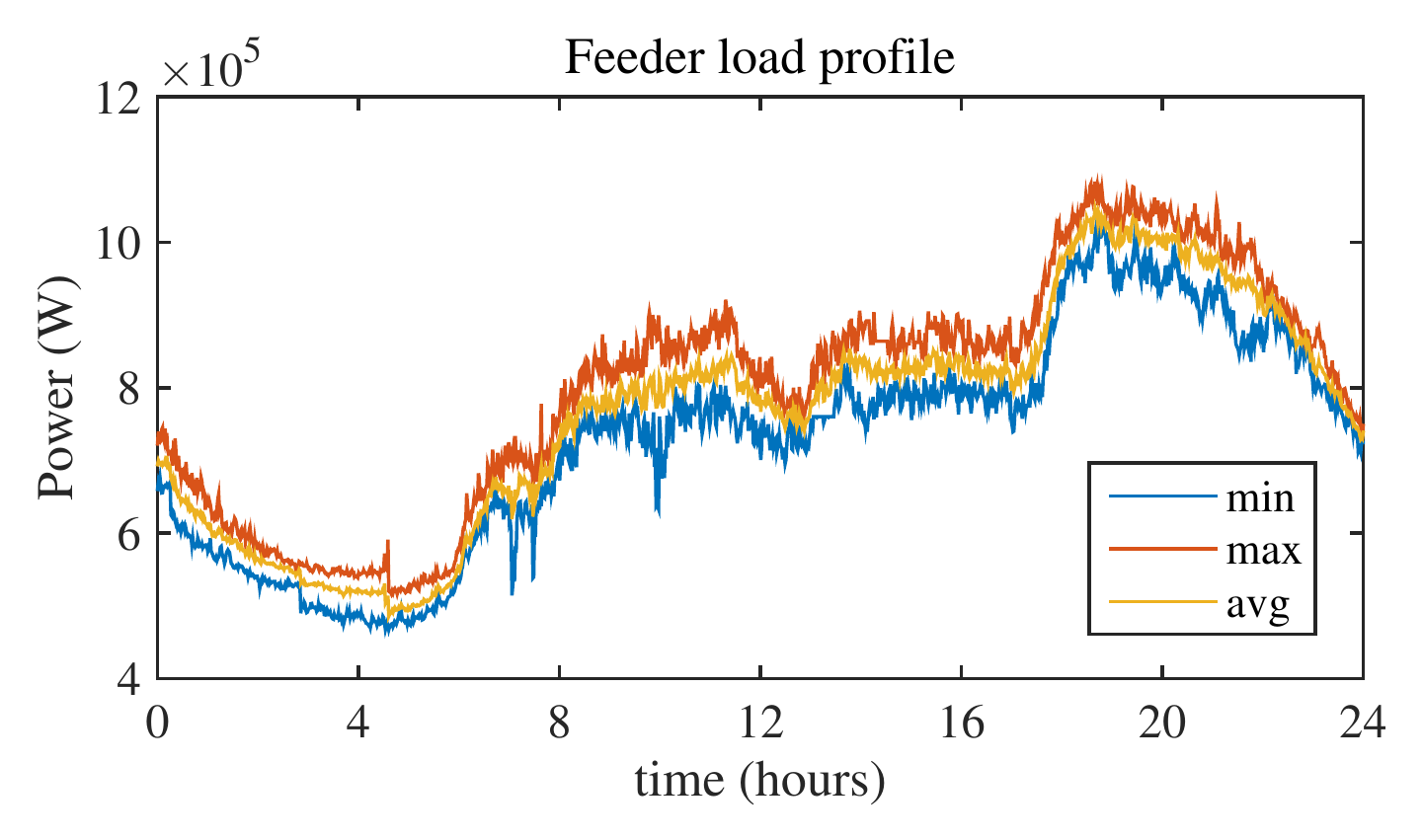}
\caption{Feeder Load Profile (Winter)}
\label{fig:3}
\end{figure}

\section{Zero Energy Reserve Approach}
\label{sec:4}

In this work, the energy storage system (ESS), generators and flexible 
loads are used to match the renewable output fluctuations. 
The measurement of solar energy $S_{\mathrm{measured}}(t)$ 
and available solar energy $S_{\mathrm{forecast}}(t+1)$ 
are formulated as presented in Section~\ref{sec:3}.

At any time, a ``difference signal'' $d(t)$ is derived from the 
$S_{\mathrm{measured}}(t)$  and  $S_{\mathrm{forecast}}(t+1)$,
which should be followed by ESS, generators and flexible loads,
\begin{align}
d(t) = S_{\mathrm{forecast}}(t+1) - S_{\mathrm{measured}}(t) 
\end{align}

A spectral decomposition of the difference signal into low and high
frequency components for the control problem on two different timescales
is proposed (not to be confused with $t$ and $t'$). The idea of spectral
decomposition is similar to the one presented for demand-side
flexibility in \cite{barooah2015}. 
The difference signal is analogous to the grid
regulation signal used to balance generation and load \cite{ela2010}. 
The high
frequency difference signal typically has a characteristic time of about 
one minute,
while the low frequency difference signal has a characteristic time 
of about five
minutes. In this work, a moving average sigmoid filter (20 samples) 
was used  to formulate the low frequency difference signal 
$d_l(t)$,
\begin{align}
d_l (t) = \sum_{k=1}^{20} w_k d(t-k)
\end{align}
where the finite impulse response was selected as a decreasing
sigmoid function $w_k = 1-(1+e^{-12 k - 0.5})^{-1}$ for 
$k \in \{1,2,\ldots,20\}$.

One can obtain a high frequency difference signal $d_h(t)$from the difference
signal and the low frequency difference signal, which can be used to
follow the minute variations in the solar output,
\begin{align}
d_h(t) = d(t) - d_l(t)
\end{align}

\begin{figure}
\centering
\includegraphics[width=1\columnwidth]{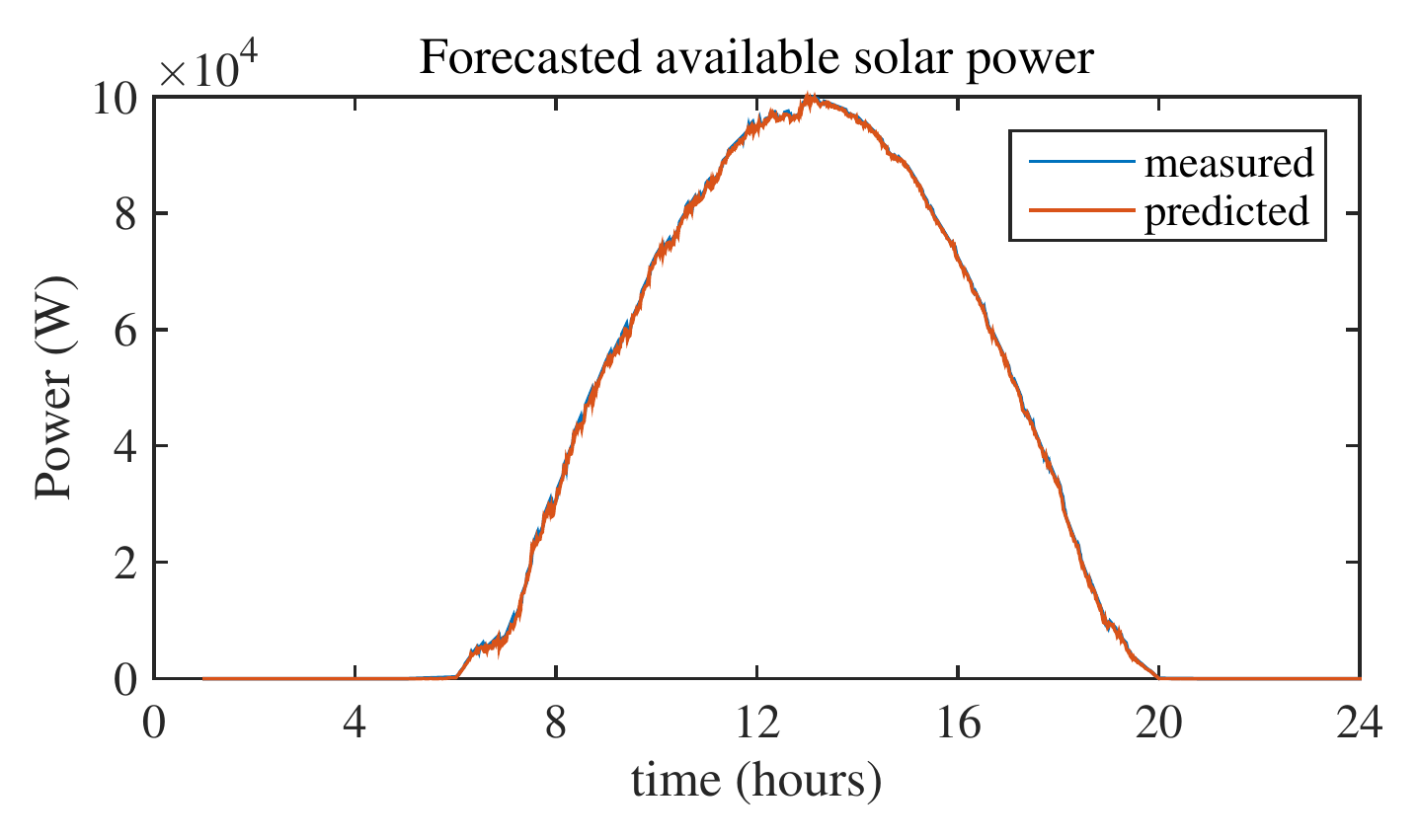}
\caption{Forecasted and Measured Solar Power}
\label{fig:4}
\end{figure}

\begin{figure}
\centering
\includegraphics[width=1\columnwidth]{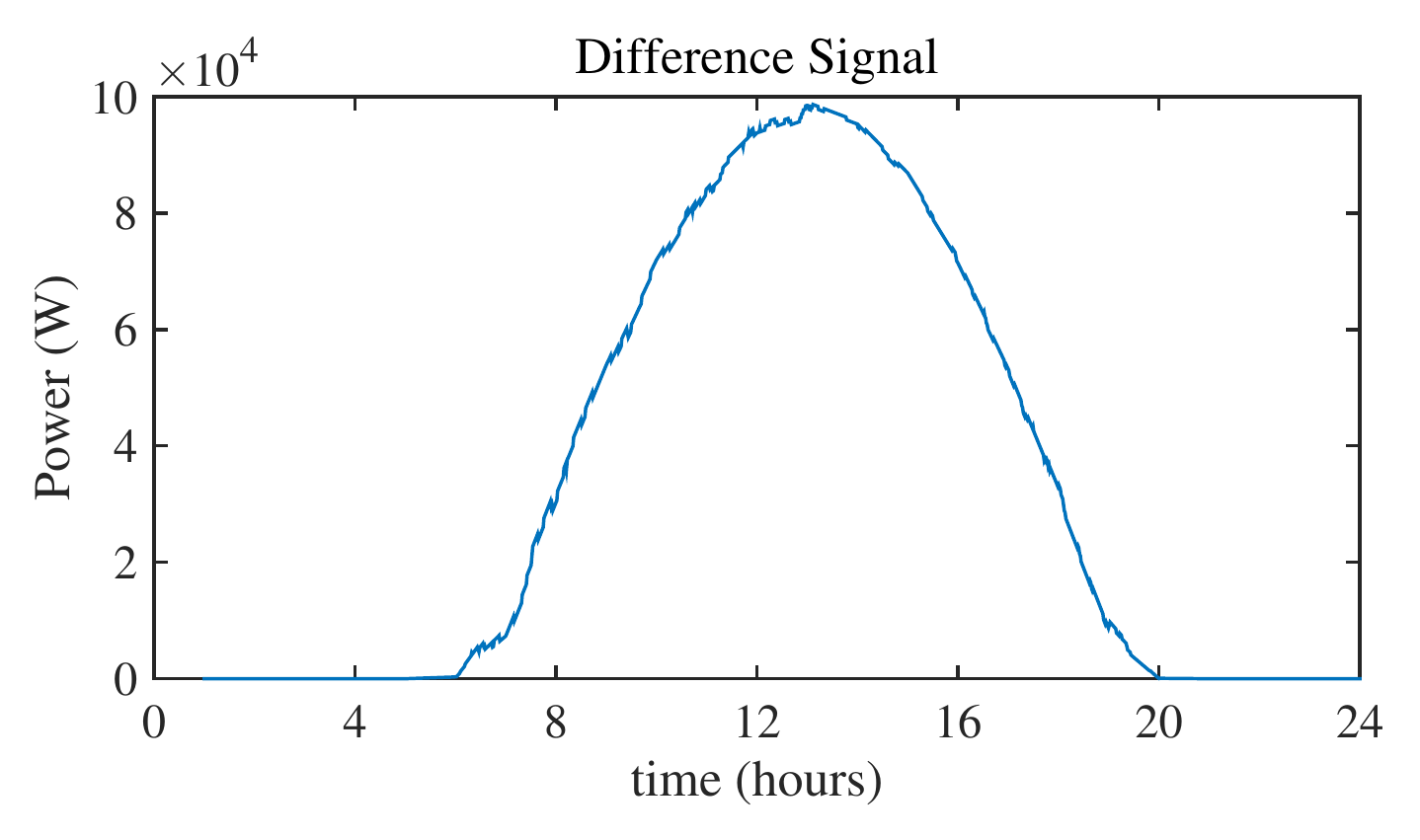}
\caption{Difference Signal}
\label{fig:5}
\end{figure}

\begin{figure}
\centering
\includegraphics[width=1\columnwidth]{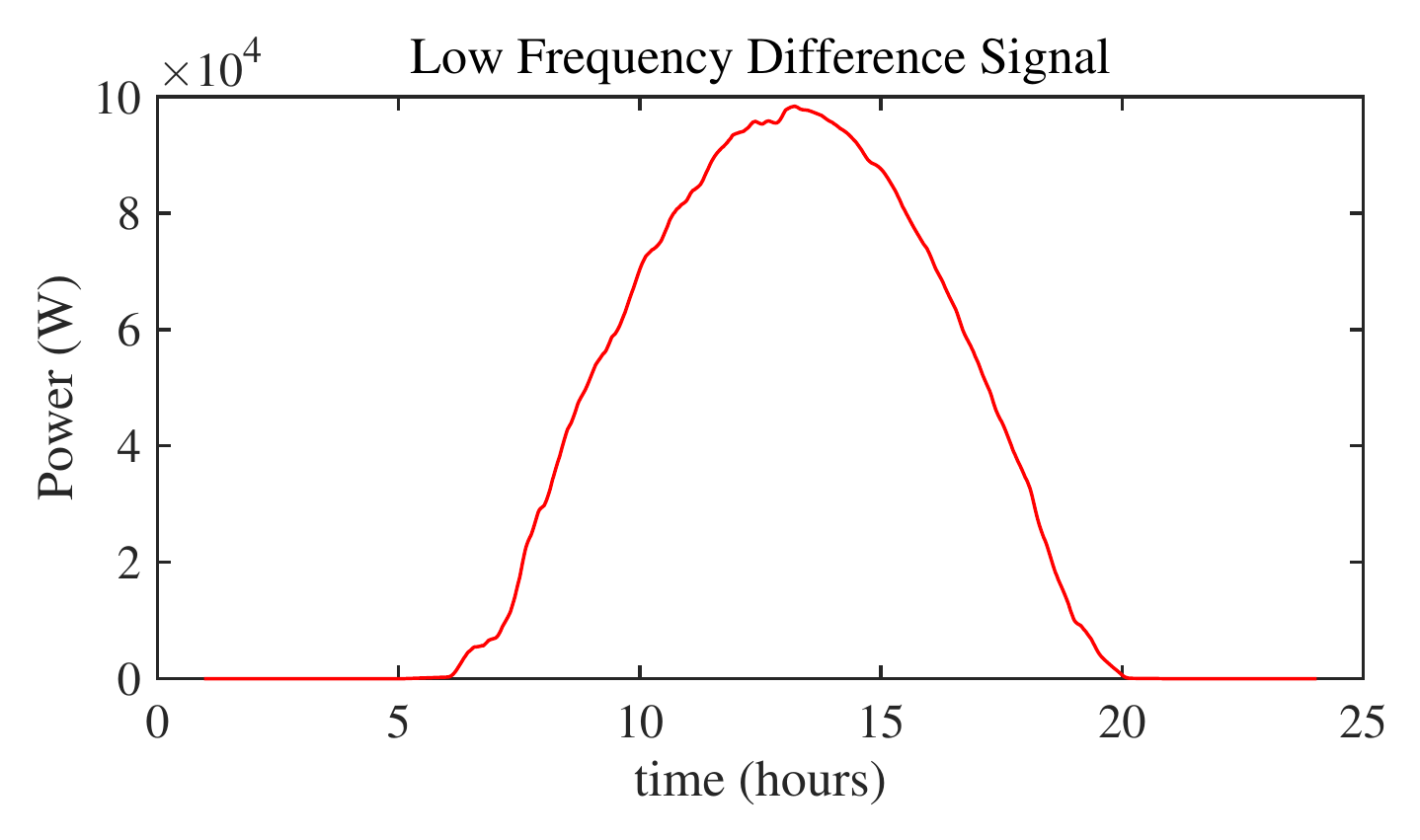}
\caption{Low Frequency Difference Signal}
\label{fig:6}
\end{figure}

\begin{figure}
\centering
\includegraphics[width=1\columnwidth]{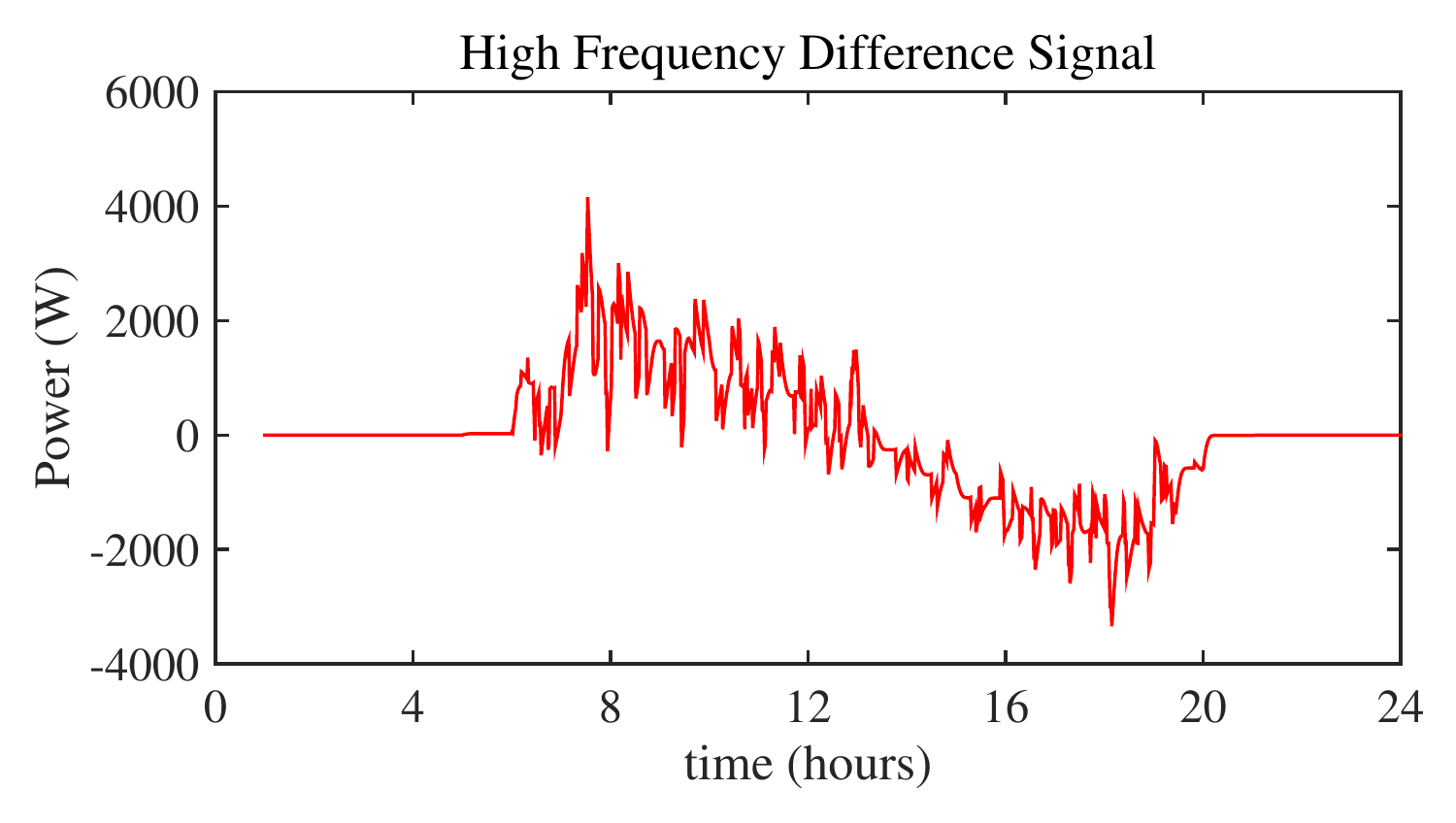}
\caption{High Frequency Difference Signal}
\label{fig:7}
\end{figure}

In Fig.~\ref{fig:4} to Fig.~\ref{fig:7} we show the forecasted and
measured solar power, the difference signal, and the low and 
high frequency difference signals.

\subsubsection{Discussion}

The averaged low frequency difference signal with a characteristic time
of about five minutes can be followed by dispatch of generators (five-minute
market) and/or modulation of power consumption of flexible load. The high
frequency difference signal can be followed by ESS, which employ
a ``zero energy reserves'' approach. The energy capacity of ESS
installed can be reduced using this control strategy because they are
only used to follow ramp ups and/or ramp downs unlike traditional
control in which ESS follows the difference signal. A signal of
frequency $r(t)/(A\pi)$ cycles/hour can be successfully tracked by the
batteries with a capacity of $A$~MWh to smooth out the intermittencies in
power caused by solar (here $r(t)$ is the maximum deviation from
base-line). The energy capacity calculations for the data used in
simulations are presented in the appendixA.

\section{Optimization Problem}
\label{sec:5}

Consider the Volt/VAR optimization problem on two timescales. As stated
in \cite{farivar2011}, 
the slow timescale problem of capacitor bank switching is to
find a state $s(t)$ of the discrete controllers (capacitor banks) which
minimized a cost function, which represents the cost of switching from
configuration $s(t-1)$ in previous time period to $s(t)$, the current time
period. This slow timescale control is used to adapt to the aggregate
reactive power requirement of the system. One can compute the optimal
setting of discrete controller using various methods studied in the
literature \cite{civanlar1985}.

$s(t)$ determines the system configuration and it is constant over $t$
while the fast timescale control time $t'$ changes.  The fast timescale
inverter control optimization problem is modeled by $C(s(t),t')$, 
which is the sum of technical losses in each phase in the distribution 
system and the sum of weighted deviations of voltage in every phase from 
nominal voltage.

\begin{align}
\min
C(s(t),t') &=
\sum_{(i,j,k) \in B \times B \times \Phi} r_{ij,k} |I_{ij,k}|^2 + 
\nonumber \\
& \qquad
\sum_{(i,k)\in\mathcal B \times \Phi} \alpha_{ik} |v_{ik}(t') - v_{\mathrm{nom}}|^2 \\
\text{subject to} \ &
|P_i + jQ_i| \leq S_i \\
& P_{\min}^L \leq P^L \leq P_{\max}^L \\
& P_{i-k}^{\mathrm{exp}} - P_{i-k} = 0 \\
& Q_{i-k}^{\mathrm{exp}} - Q_{i-k} = 0
\end{align}
where $P_{\min}^L = \min R(t)$, $P_{\max}^L = \max R^L(t)$ and
$R^L(t)$ is the variation range of the flexible load profile 
for instant $t$. The expressions for $P_{i-k}$ and $Q_{ik}$
were given in (\ref{eq:6}) and (\ref{eq:7}). Here, $\alpha_{ik}$
are the weights assigned for each bus and phase.

The term $r_{ij,k}|(I_{iji,k}|^2$ represents the
technical losses in the line between buses $i$ and $j$ in
phase $k$. In this work, ANSI C84-1, Range A (0.95 p.u to 1.05 p.u) 
is considered as the nominal voltage range. Recent studies 
\cite{schneider2010} show that reduction of energy
consumption can be achieved using Conservation Voltage Reduction (CVR).
This can be achieved by reduction of feeder voltage. In this work, the
optimization problem was solved by choosing $v_i{\mathrm{nom}}$ as 
$1$~p.u, because the goal is to achieve a flat voltage profile, thereby  
improving the $\mathrm{STARFI}_x$ PQ index \cite{dugan1996}.
The second term of optimization problem models the
average power saving which can be achieved by CVR, which is quantified
in terms of a profit function:
\begin{align}
J(t) &= K_1 \Delta C + K_2 \Delta P_l - K_3 \mathrm{STARFI}_x
\end{align}
where $\Delta C$ and $\Delta P_{l}$ are lost capacity and reduction in
technical losses of the line with 1547 control and fast inverter
control, respectively. $K_1$, $K_2$ are constants with units
\$/KW, $K_3$ is the penalty for the utility (in \$/number of customers) 
and $\mathrm{STARFI}_x$ is a long duration PQ index related to voltage
deviation \cite{dugan1996}.

\begin{figure}
\centering
\includegraphics[width=1\columnwidth]{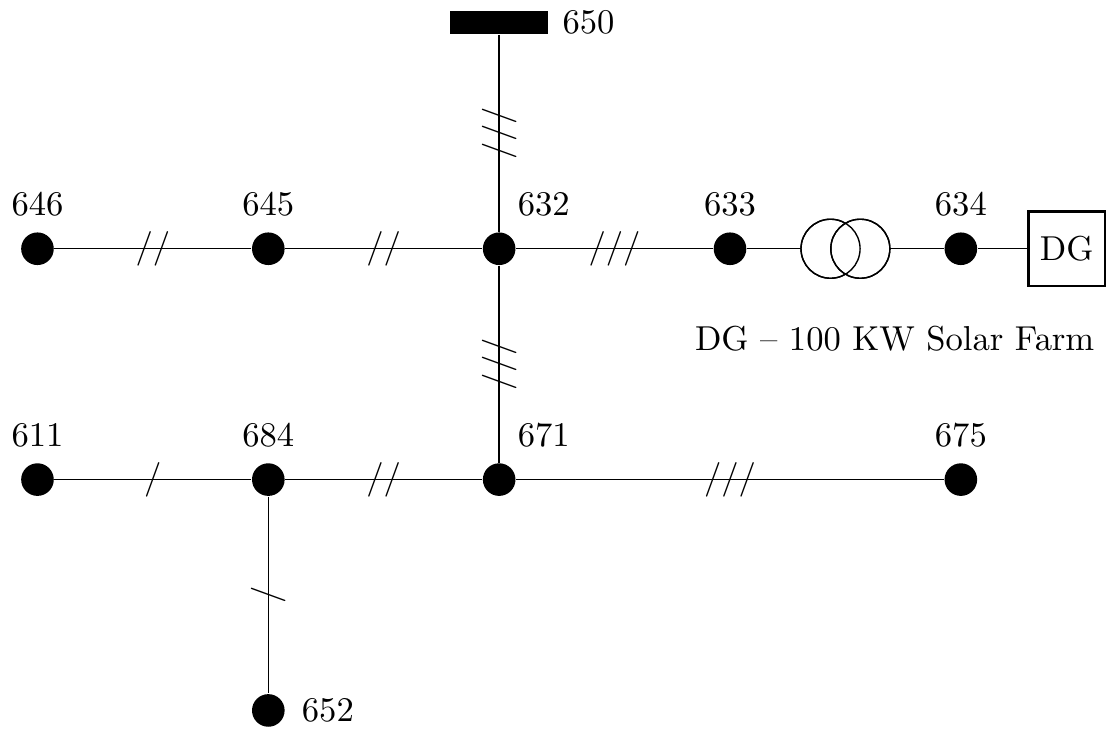}
\caption{Simplified IEEE 13 Bus System}
\label{fig:8}
\end{figure}

\begin{figure}
\centering
\includegraphics[width=1\columnwidth]{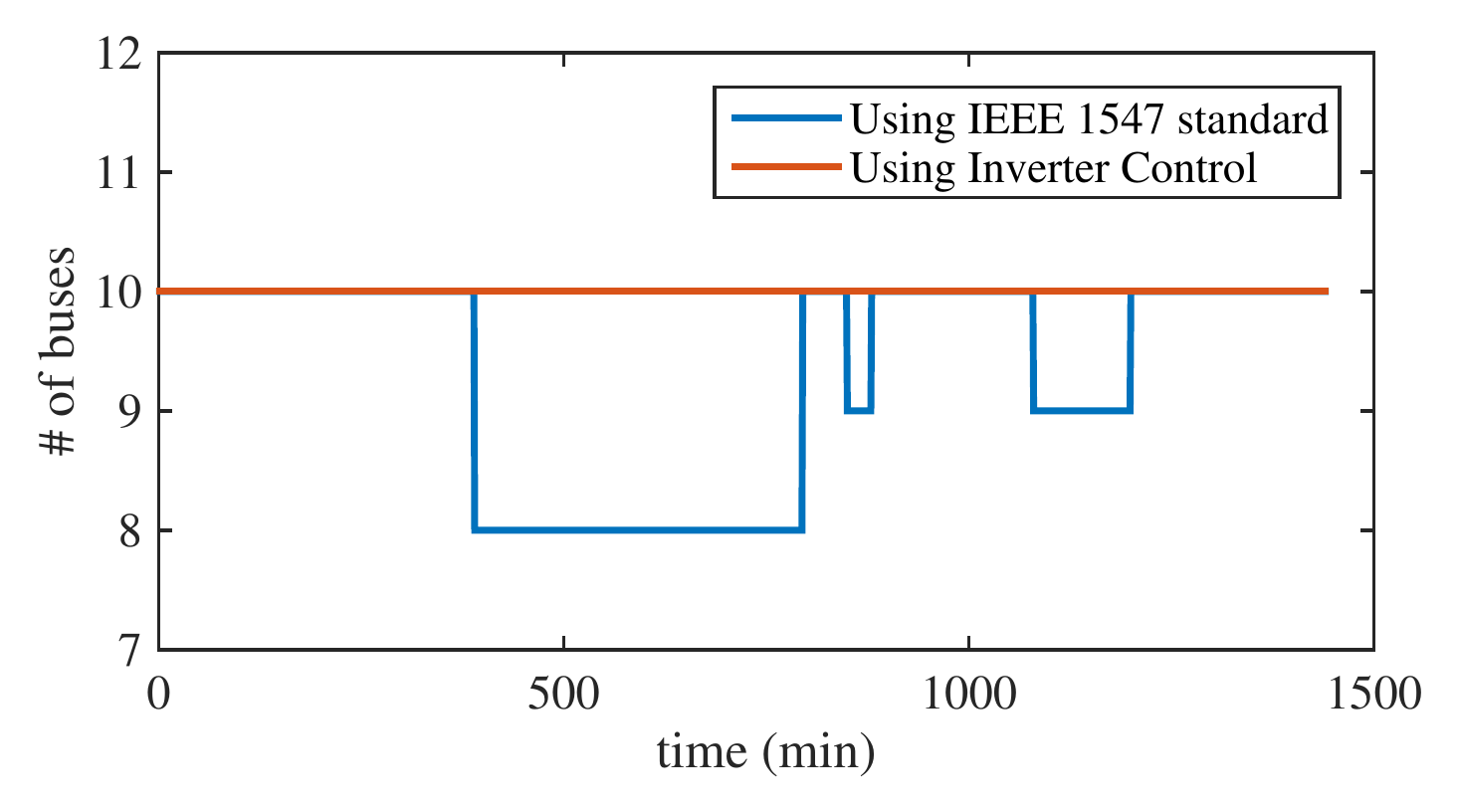}
\caption{Number of buses in the acceptable operating voltage region for
IEEE-1547 Standard and the Fast Inverter Control assuming a voltage
tolerance of 5\%.}
\label{fig:9}
\end{figure}

\begin{figure}
\centering
\includegraphics[width=1\columnwidth]{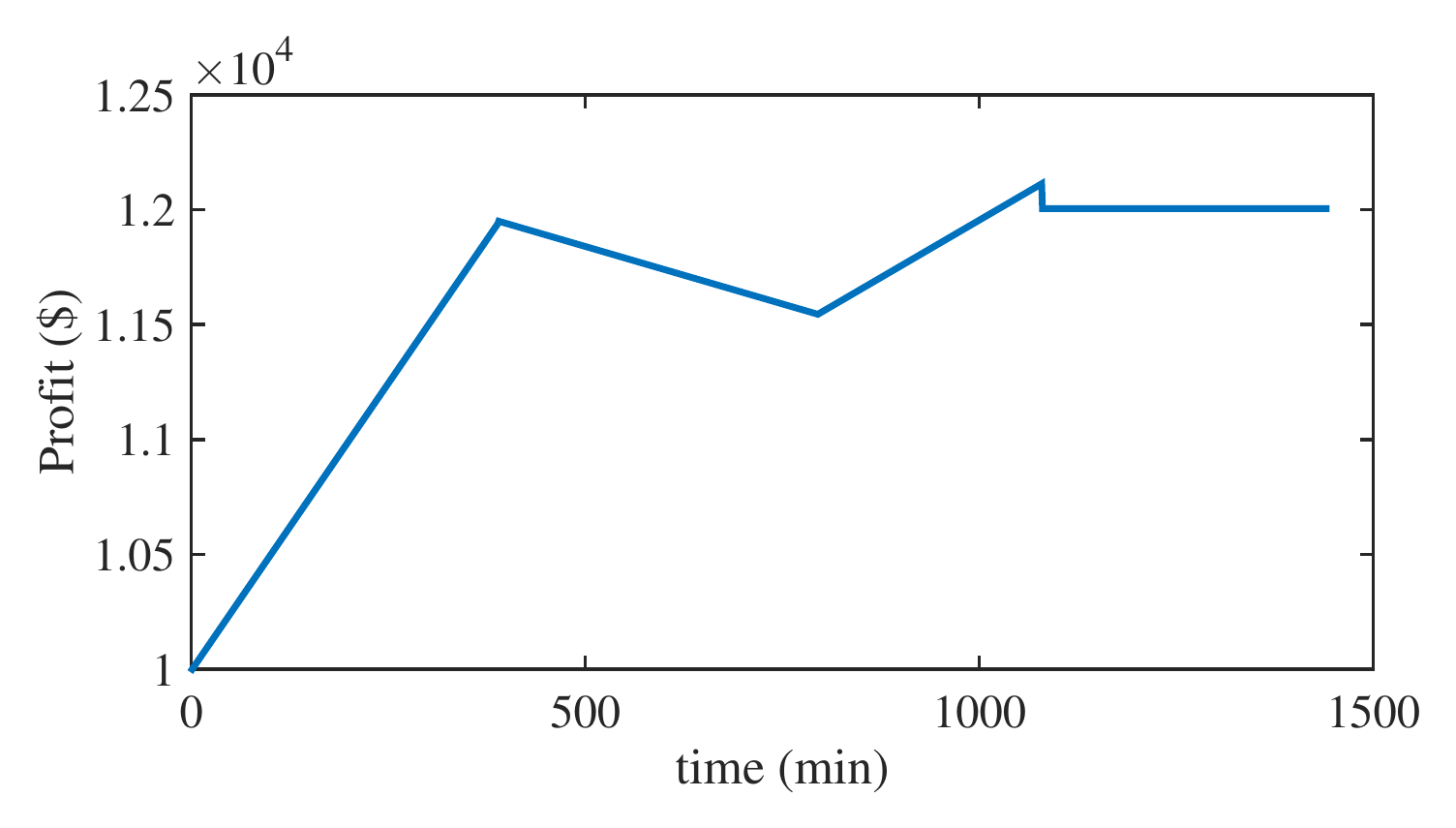}
\caption{Average Utility Profit Function}
\label{fig:10}
\end{figure}

\section{Simulation Results}
\label{sec:6}

Consider the IEEE 13 bus test feeder system, illustrated in
Fig.~\ref{fig:8}. 
 The test feeder is modified to consider renewable energy penetration
on bus 634. In this section, test  results of fast inverter control when
applied to a feeder with high renewable penetration to prove the
robustness of the control algorithm are presented.

IEEE 1547 standard states that the ``DG shall not actively regulate the
voltage and shall not cause the system voltage to go outside the
requirements of ANSI C84-1, Range A (0.95 p.u to 1.05 p.u)''. This
forbids the distributed generators to regulate the voltage. The solar
inverters operate at unity power factor w.r.t to the distribution
system. We compare our results to this current PV integration standard.

A simplified IEEE 13 Bus system as shown in Fig.~\ref{fig:8} was used. 
The bus 692 was eliminated by closing the switch and bus 680 was eliminated due
to zero power injection from the standard IEEE 13 bus system. Also, a
solar farm with $100$~{KW} of capacity was introduced at bus 634. 
The forecast of the available solar power and load forecasting was done 
as explained in Section~\ref{sec:3}.
The condition of high renewable penetration was simulated using the
appropriate profiles of solar power available as shown in
Fig.~\ref{fig:4}. The various load profiles for residential, 
industrial and commercial loads as shown in Fig.~\ref{fig:2} 
were used. These loads were distributed along the feeder. 
HVAC system in commercial buildings were used as a flexible load.
The slow timescale control was done by changing the capacitor bank
configuration to follow the reactive power requirement of the aggregated
feeder load profile.

In the following simulations, for every time step $t'$, the solar farm
was simulated as a PV bus to determine the maximum reactive power it can
inject and then simulated as a PQ bus, after the optimal operating point
has been determined by solving the optimization function as described in
Section~\ref{sec:5}.

\begin{table}
\caption{SIMULATION RESULTS FOR TRANSFORMER TAP CHANGES}
\centering
\begin{tabular}{cc}
\hline
IEEE 1547 standard & Fast Inverter VAR Control \\
43 & 19 \\
\hline
\end{tabular}
\end{table}
  
Test results demonstrate the robustness of the control strategy
presented in this paper. All the buses serving the loads were in the
acceptable range over a 24-hour period while there were some buses operating
in unacceptable ranges for significant periods of time using IEEE1547
standard.

There is significant reduction in tap changes of OLTC. They are reduced
by about half which increases the lifecycle and reduce operating costs.
The utility profit function proposed in (16), an index for savings
achieved by optimal reactive power dispatch of inverters, reducing the
system losses and peak load in terms of cost (in \$) was calculated for 
every time step $t’$. 
By assuming appropriate values for $K_1$, $K_2$ and penalty $K_3$, 
the averaged utility profit function is illustrated in Fig.~{\ref{fig:10}}.

\section{Conclusions}
\label{sec:7}
Increasing intermittent renewable energy penetration presents imperative
operation challenges to distribution utilities. In this work, a control
strategy for real-time
management of power quality design for distribution networks is
presented. The approach considers the inherent characteristics of
distribution networks such as unbalanced operation and different load
types and ESS. The presented control strategy design addresses the
induced voltage fluctuations due to variability in renewable generation.
A ``zero energy reserves'' approach to tackle fluctuations of renewable
energy generators is presented. The power consumption of flexible loads
is modulated to reduce the technical losses and peak load of the feeder
in a reactive power (VAR) control strategy. Unbalanced power flow,
different load profiles and flexible loads as “virtual energy storage”
are used to improve voltage profile and reduce system losses while
maintaining system reliability. IEEE 13 bus distribution system is used
for control strategy design validation. Comparative test results with
IEEE 1547 standard indicate reduction in the system technical losses and
the stress on automatic voltage regulators (AVR). The ease of design of
control strategy indicate potential real-life application.


\appendix

For the data used in the simulations:

\begin{flushleft}
Energy of slow difference signal $= 789.15$~KWh. \\
Energy of the fast difference signal $= 2.85\cdot 10^{-10}$~KWh. (zero energy reserve). \\
Energy Capacity of ESS to follow difference signal  $= 1.57$~KWh. \\
Energy capacity of ESS to follow fast difference signal $= 63.65$~Whr.
(to follow a maximum ramp $4.15$KW).
\end{flushleft}

\bibliography{refs}
\bibliographystyle{IEEEtran}

\end{document}